\documentclass[conference]{IEEEtran}
\IEEEoverridecommandlockouts
\usepackage[hidelinks]{hyperref}

\usepackage{cite}
\usepackage{amsmath,amssymb,amsfonts}
\usepackage{algorithm}
\usepackage{algorithmic}
\usepackage{graphicx}
\usepackage{textcomp}
\usepackage{xcolor}
\usepackage{siunitx}
\usepackage{bm}
\usepackage{tikz}
\usepackage{pgfplots}
\usepackage{comment}

\usetikzlibrary{calc}
\usetikzlibrary{fit}
\usetikzlibrary{shapes,arrows,positioning,calc,decorations.markings}
\tikzset{
  block/.style = {draw, fill=white, rectangle, minimum height=3em, minimum width=3em},
  tmp/.style  = {coordinate}, 
  sum/.style= {draw, fill=white, circle, node distance=1cm, inner sep = 0mm},
  branch/.style = {circle, fill=black,inner sep=0.35mm, outer sep = 0},
  input/.style = {coordinate},
  output/.style= {coordinate},
  pinstyle/.style = {pin edge={to-,thin,black}
  }
}
\def\BibTeX{{\rm B\kern-.05em{\sc i\kern-.025em b}\kern-.08em
    T\kern-.1667em\lower.7ex\hbox{E}\kern-.125emX}}
\pgfplotsset{compat=1.18}

\definecolor{y_color}{RGB}{60, 180, 75}
\definecolor{theta_color}{RGB}{245, 130, 48}
\definecolor{theta_dot_color}{RGB}{0, 130, 200}
\definecolor{y_dot_color}{RGB}{145, 30, 180}
\definecolor{bowling_ball_color}{RGB}{80, 207, 22}

\usetikzlibrary{angles,quotes} % for pic (angle labels)

\usepgfplotslibrary{external}
\DeclareMathOperator{\sinc}{sinc}

\usepackage{caption}
\usepackage{epsfig} % for postscript graphics files
\usepackage{booktabs}
\usepackage{colortbl}
\usepackage{multirow}
\definecolor{mygray}{gray}{0.9}
\definecolor{mydarkgray}{gray}{0.85}
\definecolor{darkgray}{RGB}{169,169,169}
\definecolor{darkgray176}{RGB}{176,176,176}
\definecolor{gray}{RGB}{128,128,128}
\definecolor{green01270}{RGB}{0,127,0}
\newtheorem{remark}{Remark}

\begin{document}

\title{Parameter Refinement of a Ballbot and Predictive Control for Reference Tracking with Linear Parameter-Varying Embedding}

\author{\IEEEauthorblockN{1\textsuperscript{st} Dimitrios S. Karachalios}
\IEEEauthorblockA{\textit{Institute for Electrical Engineering in Medicine} \\
\textit{University of Luebeck}\\
Lübeck, Germany \\
email: dimitrios.karachalios@uni-luebeck.de}
\thanks{This work was funded by the German Research Foundation (DFG),
project number 419290163.}
\and
\IEEEauthorblockN{2\textsuperscript{nd} Hossam S. Abbas}
\IEEEauthorblockA{\textit{Institute for Electrical Engineering in Medicine} \\
\textit{University of Luebeck}\\
Lübeck, Germany\\
email: h.abbas@uni-luebeck.de}
}
\maketitle
\begin{abstract}
In this study, we implement a control method for stabilizing a ballbot that simultaneously follows a reference. A ballbot is a robot balancing on a spherical wheel where the single point of contact with the ground makes it omnidirectional and highly maneuverable but with inherent instability. After introducing the scheduling parameters, we start the analysis by embedding the nonlinear dynamic model derived from first principles to a linear parameter-varying (LPV) formulation. Continuously, and as an extension of a past study, we refine the parameters of the nonlinear model that enhance significantly its accuracy. The crucial advantages of the LPV formulation are that it consists of a nonlinear predictor that can be used in model predictive control (MPC) by retaining the convexity of the quadratic optimization problem with linear constraints and further evades computational burdens that appear in other nonlinear MPC methods with only a slight loss in performance. The LPVMPC control method can be solved efficiently as a quadratic program (QP) that provides timing that supports real-time implementation. Finally, to illustrate the method, we test the control designs on a two-set point 1D non-smooth reference with sudden changes, to a 2D nonstationary smooth reference known as Lissajous curves, and to a single-set point 1D non-smooth reference where for this case theoretical guarantees such as stability and recursive feasibility are provided. 
\end{abstract}
\begin{IEEEkeywords}
Identification, Control, Ballbot, Model Predictive Control, Linear Parameter-Varying, Quadratic Program
\end{IEEEkeywords}
\section{Introduction}\label{sec:introduction}
\begin{figure}[ht!]
    \centering
    \hspace{5mm}
     \includegraphics[scale=0.03]{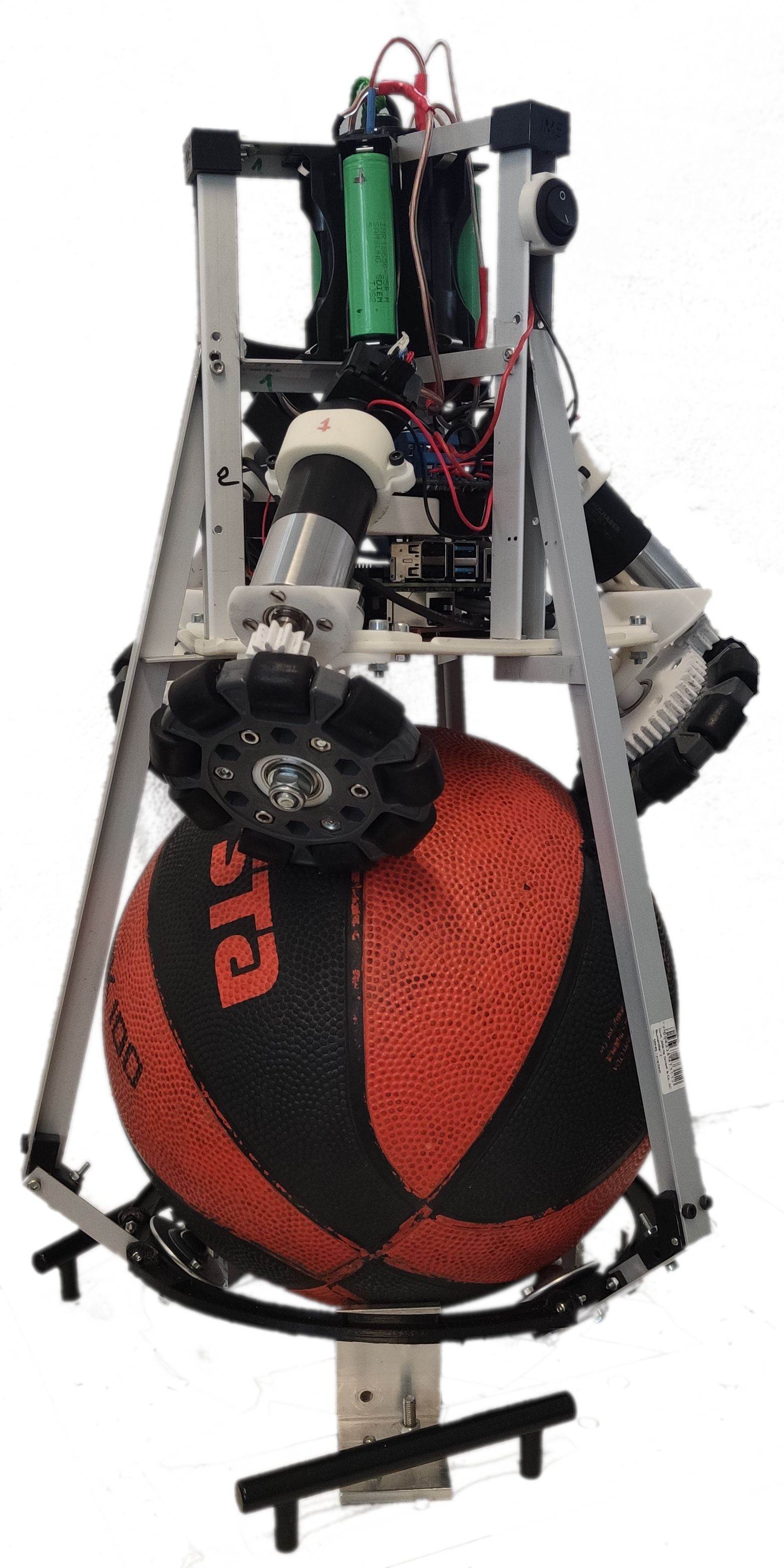} 
     \hspace{5mm}
     \includegraphics[scale=1]{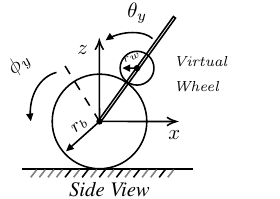}
     \vspace{2mm}
    \caption{The Ballbot constructed at the Institute for Electrical Engineering in Medicine (IME) and side view on the $xz$-plane.}
    \label{fig:ballbot_SideView}
\end{figure}

A ballbot Fig.~\ref{fig:ballbot_SideView} is an omnidirectional mobile robot that balances on a single spherical wheel \cite{Max,Peter,Pham}. The single point of contact with the ground makes this under-actuated system agile but challenging to control. The use cases of a ballbot include healthcare, virtual assistants, etc., as service robots, where maneuverability in busy environments while maintaining stability is important. A ballbot system's stability implies balancing and trajectory tracking through a predetermined path. Linear control of a ballbot \cite{Max} is proven insufficient, and the LPV model, being much closer to the nonlinear model, has shown great potential in several applications, such as autonomous vehicles \cite{Maryam}.

The dynamic model structure is well-known for robotic systems based on physical laws. The identification of the parameters can be obtained through classical mechanics \cite{SiScViOr08} and enhanced white modeling, but parameters that describe friction can be determined only experimentally. This study will provide a method that refines the parameters that define the nonlinear ballbot model from the measurements derived in \cite{Max}. For controlling the ballbot, often, PID controllers considered \cite{PuBo19,Max} or with double-loop approaches \cite{Pham}. One of the limitations of the aforementioned methods is the consideration of input/state constraints.

Model predictive control (MPC) has been widely applied for reference tracking with constraints \cite{Maryam}. It can generate optimal trajectories that steer the robot in a constrained space. Since autonomously moving robots are safety-critical systems, nonlinear MPC (NMPC) is gaining popularity due to its ability to utilize high-fidelity nonlinear models, enabling more accurate and precise control actions. Another advantage of MPC is that it can be leveraged with planning algorithms, which is crucial for such autonomously moving objects. The drawback of NMPC is the computational complexity that makes online minimization of the underlying objective function over a nonconvex manifold cumbersome. Therefore, a good alternative is to utilize the LPV formulation that can embed the nonlinear model equivalently and introduce the LPVMPC framework that can be solved efficiently as QP, with nearly the same computational burden as in linear MPC, enabling real-time implementation.

\vspace{-4mm}

\textit{Contents}: Sec.~\ref{sec:preldef} contains preliminaries and definitions. In Sec.~\ref{sec:nonlinear}, the nonlinear ballbot model is presented. Continuously, in Sec~\ref{sec:lpv}, the linear parameter-varying (LPV) model is presented with technical aspects. In Sec.~\ref{sec:stabilization}, we refer to other stabilization methods that will help tune our terminal cost, and in Sec.~\ref{sec:LPVMPC}, the LPVMPC complete framework is introduced with all the details along with algorithms for ease of implementation. Sec.~\ref{sec:Results} contains results; in particular, in Sec.\ref{sec:Identificationb}, the nonlinear model is identified. In Sec.~\ref{sec:reference_tracking}, three representative scenarios on reference tracking illustrate the developed method together with results that can ensure theoretical guarantees. Finally, a discussion followed by acknowledgment concludes the study.

\textit{Notations:} $I_n$ is the identity matrix of dimension $n$. The notation $Q\succeq 0$ represents the positive semi-definiteness of a matrix $Q$. The weighted norm $\| x\|_Q$ is defined as $\| x\|_Q^2 = x^\top Q x $ and similarly $\lVert x\rVert^{2}_P = x^\top P x,~\lVert u\rVert^{2}_R = u^\top R u$. The function $\texttt{diag}(\mathbf{x})$ constructs a diagonal square matrix from a vector $\mathbf{x}$. The set of positive integers, including zero, is denoted by $\mathbb{Z}_{+} \cup \{ 0\}$; $i|k$ denotes at time $k$ the $i^{\textrm{th}}$ prediction.

%%%%%%%%%%%%%%%%%%%%%%%%%%%%%%%%%%%%%%%%%%%%%%%%%%%%%%%%%%%%%%%%%%%
% METHODS
\section{Preliminaries and Definitions} \label{sec:preldef}
%%%%%%%%%%%%%%%%%%%%%%%%%%%%%%%%%%%%%%%%%%%%%%%%%%%%%%%%%%%%%%%%%%%
\subsection{Nonlinear Model of a Ballbot}\label{sec:nonlinear}
%%%%%%%%%%%%%%%%%%%%%%%%%%%%%%%%%%%%%%%%%%%%%%%%%%%%%%%%%%%%%%%%%%%
The ballbot model is thoroughly discussed regarding the Euler-Lagrange method of formulating the dynamics in \cite{Max,Peter}. A first attempt towards model discovery was performed in \cite{Max}, where a linearized model was used to balance the ballbot. The ballbot is a $5$-DoF system, where three are the body's rotational motion, and the two are the translation in the $xy$-plane. We will see that the body's angular position $\theta$ is used mainly for balancing, and the angular position of the ball $\phi$ is used mainly for tracking see Fig.~\ref{fig:ballbot_SideView}.

The symmetry of the ballbot control problem helps us focus independently on just the $xz$-plane \cite{Max} to simulate and understand its stability and controllability. The ballbot state space model in \eqref{eq:Euler-Lagrange} gives us the two observable and controllable states as $q = [\phi_{y}, \theta_{y}]^{T}$, together with $\tau_{y}$, being the control input torque of the virtual wheel actuator of the ballbot in the $xz$-plane as in Fig.~\ref{fig:ballbot_SideView} (right).

The ballbot model consists of the following matrices $M, C, D, G$, that denote the mass, Coriolis force (centrifugal forces), frictional torque, and gravitational forces, respectively, together with the input vector $\widetilde{B}$. The nonlinear mechanical model is 
\begin{align}
\label{eq:Euler-Lagrange}
M(q)\ddot{q} + C(q,\dot{q}) + D(\dot{q}) + G(q) = \widetilde{B}\tau_{y},
\end{align}
where the matrices are defined as:
\begin{align}
\begin{aligned}
\label{eq:parts}
    M &= \begin{bmatrix}
        b_{1} & -b_{2} + \ell r_{b}\cos{\theta_{y}} \\
        -b_{2} + \ell r_{b}\cos{\theta_{y}} & b_{3}
    \end{bmatrix}, \\
    C &= \begin{bmatrix}
        -\ell r_{b}\sin({\theta_{y})\dot{\theta}_y^{2}} \\
        0
    \end{bmatrix},~D = \begin{bmatrix}
        b_{4}\dot{\phi}_{y} \\
        0
    \end{bmatrix}, \\
    G &= \begin{bmatrix}
        0 \\
        -\ell g\sin({\theta_{y}})
     \end{bmatrix},~\widetilde{B} = \begin{bmatrix}
        \frac{r_{b}}{r_{w}} \\
        - \frac{r_{b}}{r_{w}}
    \end{bmatrix}. 
\end{aligned}
\end{align}
The nonlinear state-space representation of the planar ballbot model is given next in (\ref{eq:states}).
\begin{equation}\label{eq:states}
\begin{aligned}
    \underbrace{\begin{bmatrix}
        \dot{q} \\
        \ddot{q}
    \end{bmatrix}}_{\dot{x}}
    =\underbrace{\begin{bmatrix}
        \dot{q} \\
        M^{-1}(-C-D-G-\epsilon) 
    \end{bmatrix} + \begin{bmatrix}
        0 \\
        M^{-1}\widetilde{B}
    \end{bmatrix}\tau_y}_{f(x,u)},
\end{aligned}
\end{equation}
 where, $x = [\phi_y, \theta_y, \dot{\phi}_y, \dot{\theta}_y]^\top$ is the state vector of the system with $\phi_y$ the angular displacement, $\theta_y$ the tilt angle, $\dot{\phi}_y$ the angular speed of the ballbot and $\dot{\theta}_y$ the angular speed of the tilt angle. The input to each plane, e.g., $\tau_y$, is the virtual wheel's torque Fig.~\ref{fig:ballbot_SideView}. The virtual torques $(\tau_x,\tau_y,\tau_z)$ can be converted to real torques of the three motors $(\tau_1,\tau_2,\tau_3)$ as in \cite{Max} by using \eqref{eq:torque}, where for the design in Fig.~\ref{fig:ballbot_SideView} it holds that the Zenith angle is $\alpha=\pi/4=45^{\circ}$.
 \begin{equation} 
 \label{eq:torque}
     \begin{bmatrix}
         \tau_{1} \\
         \tau_{2} \\
         \tau_{3}     \end{bmatrix} = \frac{1}{3}\begin{bmatrix}
         \frac{2}{\cos{\alpha}} & 0 & \frac{1}{\sin{\alpha}} \\
         -\frac{1}{\cos{\alpha}} & \frac{\sqrt{3}}{\cos{\alpha}} & \frac{1}{\sin{\alpha}} \\
         -\frac{1}{\cos{\alpha}} & -\frac{\sqrt{3}}{\cos{\alpha}} & \frac{1}{\sin{\alpha}} 
         \end{bmatrix}\begin{bmatrix}
         \tau_{x} \\
         \tau_{y} \\
         \tau_{z}     \end{bmatrix}.
 \end{equation}
Now that we have a well-defined nonlinear model, we can introduce an equivalent representation with an LPV embedding.
\subsection{The Linear Parameter-Varying (LPV) Embedding}\label{sec:lpv}
The linear parameter-varying (LPV) embedding \cite{Cisneros_2016,MORATO202064,Abbas2016,Abbas24} can be seen as a structured linearization over an adaptive operational point that depends on the scheduling parameters. The nonlinear system in \eqref{eq:states}, can be embedded in the continuous-time LPV representation by denoting with the subscript "c" the continuous operators that depend on the variable $\rho$ as $A_c(\rho)$, $B_c(\rho)$ and define the continuous-time LPV model as

\begin{equation}\label{eq:lpv}
\left\{\begin{aligned}
    \dot{x}(t) &= A_{c}(\rho(t))x(t) + B_{c}(\rho(t))u(t) \equiv f(x(t),u(t)),\\
    \rho(t) &= \sigma(x(t)),~x_{0} = x(0),~t \geq 0,
\end{aligned}\right.
\end{equation}
where $\rho$ is the so-called scheduling parameter that belongs to the scheduling set $\mathcal{P}\subset\mathbb{R}^{n_\rho}$ that defines the range of the scheduling parameter, and it should be compact. Further, $\rho$ depends on the states through the mapping $\sigma\colon\mathbb{R}^{n_{x}} \rightarrow \mathbb{R}^{n_{\rho}}$. 
% \begin{equation}
%     \sigma(x)=\left[\begin{array}{cccc}
%     0 & 1 & 0 & 0 \\
%     0 & 0 & 0 & 1
% \end{array}\right]\left[\begin{array}{c}
%      \phi  \\
%      \theta \\
%      \dot{\phi} \\
%      \dot{\theta}
% \end{array}\right]=\left[\begin{array}{cc}
%    \theta  & \dot{\theta} \\
% \end{array}\right]
% \end{equation} 

For the nonlinear system under consideration (on one plane of the ballbot), the state dimension is $n_x=4$. We consider two scheduling parameters of the LPV model as $\theta_y(t)$ and $\dot{\theta}_y(t)$, thus $n_\rho = 2$ and $\rho(t)=(\theta_y(t),\dot{\theta}_y(t))$. Therefore, in this case, the functionality $\sigma$\footnote{Functionality: $\sigma(x)=x^\top\left[\begin{array}{cccc}
    0 & 1 & 0 & 0 \\
    0 & 0 & 0 & 1
\end{array}\right]^\top=\left[\begin{array}{cc}
   \theta_y  & \dot{\theta}_y \\
\end{array}\right]=\rho$.} can be seen as a linear function that chooses from the state vector only some states. Finally, the input dimension is $n_u=1$ with $u(t)=\tau_y(t)$ representing the virtual torque of the virtual wheel Fig.~\ref{fig:ballbot_SideView}.

The considered LPV system's matrices of the continuous-time state-space model are given by \eqref{eq:LPV_cont}, and the elements of the matrices $A_c$ and $B_c$ are given in Table \ref{tab:elements}.

\begin{equation}
\label{eq:LPV_cont}
    \begin{aligned}
    \footnotesize
        A_{c}(\rho)=\begin{bmatrix}
        0 & 0 & 1 & 0 \\
        0 & 0 & 0 & 1 \\
        0& A_{32} & A_{33} & A_{34}  \\
        0 & A_{42}  & A_{43} & A_{44} 
        \end{bmatrix},~B_{c}(\rho)=\begin{bmatrix}
        0 \\
        0 \\
        B_{31} \\
        B_{41}
    \end{bmatrix},
    \end{aligned}
\end{equation}
\begin{table}[!h]
\centering
\caption{Elements of the operators $A_c(\rho)$ and $B_c(\rho)$, with their respective scheduling parameters given within the parentheses $A_{ij}(\rho)$ and $B_{ij}(\rho)$, and their scheduling signals are given by $\rho=(\theta_y, \dot{\theta}_y)$ as shown in the table.}
\vspace{4mm}
\begin{tabular}{|c|c|} 
\hline\\[-3mm]
 %$\rho = \sigma\big(\theta_y, \dot{\theta}_y\big)$ & \textbf{Scheduling signals}\\[2mm]
 % \hline & \\[-2mm]
  $d(\theta_{y})$ & $b_{1}b_{3} - (-b_2+\ell r_{b}\cos({\theta_{y}}))^{2}$\\[2mm]
  \hline & \\[-2mm] $A_{32}\big(\theta_y, \dot{\theta}_y\big)$ & $\frac{\ell g\sinc(\theta_{y})(b_{2}-\ell r_{b}\cos({\theta_{y}}))}{d(\theta_{y})}$\\[2mm]
  \hline & \\[-2mm]
  $A_{33}\big(\theta_y, \dot{\theta}_y\big)$ & $-\frac{b_{3}b_{4}}{d(\theta_{y})}$ \\[2mm]
  \hline & \\[-2mm]
  $A_{34}\big(\theta_y, \dot{\theta}_y\big)$ & $\frac{b_{3}\ell r_{b}\sin({\theta_{y}})\dot{\theta}_{y}}{d(\theta_{y})}$ \\[2mm]
  \hline & \\[-2mm]
 $A_{42}\big(\theta_y, \dot{\theta}_y\big)$ & $\frac{b_{1}\ell g\sinc(\theta_{y})}{d(\theta_{y})}$ \\[2mm]
  \hline & \\[-2mm]
 $A_{43}\big(\theta_y, \dot{\theta}_y\big)$ & $-\frac{b_{4}(b_{2}-\ell r_{b}\cos{(\theta_{y})})}{d(\theta_{y})}$ \\[2mm]
  \hline & \\[-2mm]
  $A_{44}\big(\theta_y, \dot{\theta}_y\big)$ & $\frac{\ell r_{b}\sin{(\theta_{y})}\dot{\theta}_{y}(b_{2} - \ell r_{b}\cos{(\theta_{y})})}{d(\theta_{y})}$ \\[2mm]
  \hline & \\[-2mm]
 $B_{31}\big(\theta_y, \dot{\theta}_y\big)$ & $\frac{r_{b}(b_{3}-(b_{2}-\ell r_{b}\cos{(\theta_{y})})}{r_{w}d(\theta_{y})}$\\[2mm]
  \hline & \\[-2mm]
 $B_{41}\big(\theta_y, \dot{\theta}_y\big)$ & $\frac{r_{b}((b_{2}-\ell 
 r_{b}\cos{(\theta_{y})})-b_{1})}{r_{w}d(\theta_{y})}$\\[2mm]
  \hline 
\end{tabular}
\label{tab:elements}
\end{table}
A discretization of the nonlinear continuous model is mandatory before setting the framework that combines the LPV with Model Predictive Control (MPC). Therefore, to discretize \eqref{eq:lpv}, many methods exist that differ in accuracy and numerical properties (s.a., numerical stability). For instance, the Forward Euler method is commonly used due to its ease of implementation, but it usually endures low-level accuracy and, most importantly, inherits numerical instability. Therefore, we evade such numerical pitfalls by utilizing the Runge–Kutta $4^\textrm{th}$ order with a single-time step (RK4). Consider the sampling time $t_s$ that spans the continuous time as $t=t_sk,~k\in\mathbb{Z}_+\cup\{0\}$, and, it is denoted as $x(t_k) = x(kt_{s}) = x_{k}$. We denote the angular velocity of the tilt angle $\theta_y$ as $\dot{\theta}_y(t):=\omega_y(t)$ and the corresponding discrete value is $\omega_k$. Similarly, we denote the angular velocity $\dot{\phi}_y(t):=\varphi_y(t)$, thus, the complete discrete state vector remains $x_k=\left[\begin{array}{cccc}
    \phi_k & \theta_k & \varphi_k & \omega_k\end{array}\right]^\top$.

By denoting the continuous in-time LPV operator $f_{\text{LPV}}$ (i.e., $f_{\text{LPV}}$ is equivalent with the nonlinear operator $f$) of the ballbot system, we can rewrite \eqref{eq:lpv} as $\dot{x}(t)=f_{\text{LPV}}(x(t),u(t)):=A_c(\sigma(x(t)))x(t)+B_c(\sigma(x(t)))u(t)$ with the input $u(t)=\tau_y(t)$. The Runge-Kutta discretization method of order $4^{\textrm{th}}$ (RK4) with a single time step $t_s$ when zero-order hold (ZOH) is considered remains:
\begin{equation}\label{eq:dlpv}
\text{RK4}~\left\{\begin{aligned}
         \kappa_1&=f_{\text{LPV}}\left(x_k,u_k\right)\\
         \kappa_2&=f_{\text{LPV}}\left(x_k+(t_s/2)\kappa_1,u_k\right)\\
         \kappa_3&=f_{\text{LPV}}\left(x_k+(t_s/2)\kappa_2,u_k\right)\\
         \kappa_4&=f_{\text{LPV}}\left(x_k+t_s\kappa_3,u_k\right)\\
    x_{k+1} &= x_{k}+(t_s/6)(\kappa_1+2\kappa_2+2\kappa_3+\kappa_4). \\
    \end{aligned}\right.
\end{equation}
 Thus, we introduce the following remark to proceed with linear operations with the LPV formulation, which will be essential for deriving a predictor that retains the characteristics of the underlying optimization problem.
\begin{remark}[The operator $f_{\text{LPV}}$ with a given scheduling variable $\rho$ is linear]\label{rem:LPVoperator}
When the scheduling parameter $\rho$ is given, the LPV operator remains linear w.r.t. the state $x$ and input $u$. In such cases, the LPV continuous operator will be denoted as:
\begin{equation}
    f_{\text{LPV}}^{(\rho)}(\rho,x(t),u(t)):=A_c(\rho)x(t)+B_c(\rho)u(t).
\end{equation}
\end{remark}
Consequently, when the scheduling parameter can be considered known, we will denote the LPV operator with a superscript $\rho$ as $f_{\textrm{LPV}}^{(\rho)}$ that will indicate that $\rho$ is independent of the state/input allowing linear actions on both. In particular, Remark~\ref{rem:LPVoperator} allows the LPV predictor to be inserted in a linear MPC framework and retain all the characteristics of the optimization problem that stay invariant under linear operations. For instance, the convexity of the optimization problem and the linear constraints can be retained. At the same time, efficient algorithms of QP can solve the underlying control problem efficiently. 
\subsection{Stabilization and closed-loop parameters identification}\label{sec:stabilization}
To identify the parameters of the ballbot as a real physical system, it was stabilized in \cite{Max}  through a proportional integral derivative (PID) controller. During stabilization in \cite{Max}, a multi-harmonic excitation signal provided enough measurements to identify a linear model that is defined with the $p$-parameters provided in Table \ref{tab:pparameters} and in the study \cite{Max}. The identified linear model could stabilize the ballbot with an optimal feedback gain $K$ utilizing the solution of the linear quadratic regulator (LQR) approach.
\begin{table}[!h]
\centering
\caption{Physical parameters of the ballbot model  (\ref{eq:parts}), with $p$ being the linearized model parameters from \cite{Max}.}
\footnotesize
\begin{tabular}{|c |c |c |c | c |}
\hline 
$p_1$ & $p_2$ & $p_3$ & $p_4$ & $p_5$  \\ [1ex]  
-342.6038  &  -52.8301 & -1425.9 & -36.0734 & -9.1477  \\ [1mm] 
 \hline
   $p_6$ & $\ell[m]$&  $r_b[m]$ & $r_w[m]$ & $g[m/s^2]$   \\ [1ex] 
   251.8476 & 0.2978 & 0.12 & 0.05 & 9.81   \\
    \hline 
\end{tabular}
\label{tab:pparameters}
\end{table}
In this study, we will use the identified linear model from \cite{Max} for discovering the $b$ parameters that define the refined values in \eqref{eq:Euler-Lagrange} after solving (offline) a nonlinear optimization problem as will be explained in Sec.~\ref{sec:Identificationb}. Consequently, the LPV embedding of the identified nonlinear model will allow the LPV integration within the proposed MPC scheme as shown in Sec.~\ref{sec:reference_tracking} that will handle the stabilization and reference tracking of the ballbot simultaneously without cascading complex control designs based on the linearized model.
\subsection{Model Predictive Control with LPV Embeddings}\label{sec:LPVMPC}
Model predictive control with linear parameter-varying embeddings implies that the  MPC optimally controls an underlying time-varying model based on the scheduling parameters shown in the LPV formulation \eqref{eq:lpv}. 
The standard MPC form begins with the cost function defined as $J_{k}(u)$ as seen in (\ref{eq:energy}), where $x_{i|k}$ denotes the discretized state vector at time $k$ and prediction $i$. The $i$ varies in length $i=0,\ldots,N-1$ with $N$ the prediction horizon length. The reference state vector is given as $x^{\mathrm{ref}}_{i|k}$. 

\begin{equation}
    \label{eq:energy}
    \begin{aligned}
    J_k(u) &= \sum_{i=0}^{N-1}(\lVert x_{i|k}- x^{\mathrm{ref}}_{i|k}\rVert^{2}_Q + \lVert u_{i|k}\rVert^{2}_R)+\\
    &~~~~~~~~~~+\underbrace{(\lVert x_{N|k}- x^{\mathrm{ref}}_{N|k}\rVert^{2}_P}_{\text{terminal cost}}, 
    \end{aligned}
\end{equation}
\noindent
where $Q\succeq 0,~R\succ 0,~P\succ 0$ are weighting symmetric matrices. Theoretically, to provide closed-loop stabilization of the LPV model, for some feedback gain $K(\rho)$ the $A(\rho)-B(\rho)K(\rho)$ should be Hurwitz for all possible parameter variation of $\rho$. In this study, with the suitable tuning of the LQR on the linearized model, we could provide a non-parametric state feedback gain $K$ that satisfies the Hurwitz condition for the LPV over a dense grid that covers the scheduling space $\mathcal{P}$. In addition, the algebraic Riccati solution $P$ will be used to penalize the terminal cost of the MPC.

In \eqref{eq:energy} and on the one hand, the explicit decision variables are the input entries $u_{i|k},~i=0,\ldots,N-1$. Thus, the optimizer should provide a solution in $u$ that minimizes the energy function \eqref{eq:energy}. On the other hand, the $x_{i|k},~i=0,\ldots,N$ consists of the set with the implicit decision variables inserted through a predictor. Suppose the prediction $x_{i|k}$ depends nonlinearly on previous states. In that case, the energy function will not remain quadratic or convex, inevitably leading to nonconvex optimization formulation where the cost of solving such problems inherits high complexity. 

To evade such a computational burden, by utilizing the LPV formulation together with the property in Remark~\ref{rem:LPVoperator}, when the scheduling variable is given, the implicit decision variables $x_{i|k}$ are introduced linearly in the energy function and retain the quadratic manifold where the minimization is quite efficient through the quadratic program (QP) that will introduce next. 

Importantly, good scheduling estimation will make the adaptive nature of the LPV over a potentially nonlinear manifold quite accurate. Good predictions for the scheduling variable lie at the heart of the LPV method, and methods that introduce ways, e.g., scheduling tube control \cite{Abbas24}, can provide further guarantees, such as recursive feasibility and stability.    

Back to our aim, we want to minimize the quadratic cost function $J$ in \eqref{eq:energy} under the linear equality constraints that satisfy the LPV (given scheduling) predictor model together with the following linear inequalities given in \eqref{eq:constraints}.
\begin{align}
    \begin{aligned}
    \label{eq:constraints}
      &\mathcal{X} = \{ x_k \in \mathbb{R}^{n_{{\mathrm{x}}}}~|~G^x x_k \leq h^x \},\\
       &\mathcal{U} = \{ u_k \in \mathbb{R}^{n_u}~|~G^u u_k \leq h^u \},
    \end{aligned} 
\end{align}
where, the state/input inequality constraints emerge as
\begin{equation*}
    \begin{aligned}
        G^x&=\left[\begin{array}{c}
         +I_{n_{\mathrm{x}}} \\
         -I_{n_{\mathrm{x}}}
      \end{array}\right]\in\mathbb{R}^{(2n_{\mathrm{x}})\times n_{\mathrm{x}}},~h^{x}=\left[\begin{array}{c}
           +x_{\text{max}}  \\
           -x_{\text{min}} 
      \end{array}\right]\in\mathbb{R}^{2n_\mathrm{x}},\\
      G^u&=\left[\begin{array}{c}
         +I_{n_{u}}\\
         -I_{n_{u}}
      \end{array}\right]\in\mathbb{R}^{(2n_u)\times n_{u}},~h^{u}=\left[\begin{array}{c}
           +u_{\text{max}}  \\
           -u_{\text{min}} 
      \end{array}\right]\in\mathbb{R}^{2n_u}.
    \end{aligned}
\end{equation*}
The superiority of the MPC control strategy comes exactly at the point where constraints can be introduced. MPC casts the control problem into optimization with the major feature, among other methods, that can handle constraints. MPC can be extended easily to include trajectory planning of the ballbot when it moves in a structured environment with obstacles. Tangent planes can play a significant role in keeping the optimization problem as QP, like the work in \cite{Maryam}. 

The QP optimization problem is formulated together with the LPV predictor at the "present" time $k$ as the function $\mathrm{QP}(\rho_{i|k}, x_k, x^{\mathrm{ref}}_{i+k})$ with the input arguments; 1) the given (from previous time) scheduling signal $\rho_{i|k}$; 2) the initial conditions $x_k$; 3) the reference $x_{i+k}^{\textrm{ref}}$. The Quadratic Program's (QP) functionality is defined next in Algorithm~\ref{algo:LPVQPOptimProb}.

\begin{algorithm}
    \caption{The quadratic program (QP) based-LPVRK4}
\begin{algorithmic}[1]
 \renewcommand{\algorithmicrequire}{\textbf{Input:}}
 \renewcommand{\algorithmicensure}{\textbf{Output:}}
 \REQUIRE The time step $k\in\mathbb{Z}_+$, the weighted costs $Q,~R,~P$, the triplet $(\rho_{i|k}, x_k, x^{\mathrm{ref}}_{k+i})$, for $i=0,\ldots,N-1$, the horizon length $N\in\mathbb{Z}_+$ and the sampling time $t_s\in\mathbb{R}_+$.
 \ENSURE Optimal numerically design $u_{i|k},~i=0,\ldots,N-1$.\\
 \textit{Minimization of the quadratic cost function:}\\[2mm]
\hspace{-5mm}$\underset{u_{i|k}}{min} \sum_{i=0}^{N-1}\Big(\lVert\hat{x}_{i|k}-x^{\mathrm{ref}}_{k+i}\rVert^{2}_Q+\lVert u_{i|k}\rVert^{2}_R\Big)+\underbrace{\lVert \hat{x}_{N}-x^{\mathrm{ref}}_{k+N}\rVert_P^2}_{\text{terminal cost}}$\\
\textit{Subject to:}
\STATE $\hat{x}_{0|k}=x_{0|k}=x_k$,
\FOR{i=0,\ldots,N-1}
\STATE $\kappa_1=f_{\textrm{LPV}}^{(\rho_{i|k})}\left(\rho_{i|k},\hat{x}_{i|k},u_{i|k}\right)$,\\
\STATE  $\kappa_2=f_{\textrm{LPV}}^{(\rho_{i|k})}\left(\rho_{i|k},\hat{x}_{i|k}+(t_s/2)\kappa_1,u_{i|k}\right)$,\\
\STATE  $\kappa_3=f_{\textrm{LPV}}^{(\rho_{i|k})}\left(\rho_{i|k},\hat{x}_{i|k}+(t_s/2)\kappa_2,u_{i|k}\right)$,\\
\STATE  $\kappa_4=f_{\textrm{LPV}}^{(\rho_{i|k})}\left(\rho_{i|k},\hat{x}_{i|k}+t_s\kappa_3,u_{i|k}\right)$,\\
\STATE Update: $\hat{x}_{i+1|k} = \hat{x}_{i|k}+(t_s/6)(\kappa_1+2\kappa_2+2\kappa_3+\kappa_4)$,\\
\ENDFOR
\STATE Satisfy: $\hat{x}_{i|k} \in \mathcal{X}, \forall i = 1,\dots,N$,
\STATE Satisfy: $u_{i|k} \in \mathcal{U}, \forall i = 0,1,\dots,N-1$.
\end{algorithmic}
\label{algo:LPVQPOptimProb}
\end{algorithm}
%%%% Explanation of the algorithm 1
Next, we discuss Algorithm~\eqref{algo:LPVQPOptimProb} step by step when it is called at the simulation time instance $k\in\mathbb{Z}_+$:
\begin{itemize}
    \item At time $k$, the input to the QP consists of the matrices $Q,~R,~P$, the initial condition $x_k$, the reference $x_{k+i}^{\textrm{ref}}$, and the scheduling variable that is either available from previous time $k-1$ or initialized from initial conditions $x_k$. In particular, when the simulation starts at time $k=0$, and there is no available scheduling prediction from the past, the scheduling variable is initialized as $\rho_{i|0}=\sigma(x_0),~i=0,\ldots,N-1$.
    \item The optimization problem consists of the quadratic cost function \eqref{eq:energy} that is updated with the implicit decision variables $\hat{x}_{i|k}$ as a function of the explicit decision variables $u_{i|k}$ with the LPV predictor. This can be done with linear operations since the scheduling variable $\rho_{i|k},~i=0,\ldots,N-1$ is considered known, and the linear property explained in Remark \ref{rem:LPVoperator}. The RK4 in lines (2-7) is the discretization scheme that predicts the evolution.
    \item The QP can be solved together with the constraints in lines (9-10), efficiently from \cite{Lofberg2004} and provide the optimized input design at each time $k$ numerically as $u_{i|k},~i=0,\ldots,N-1$.
\end{itemize}
 %%%%%%%%%%%%%%%%%%%%%%%%%%%%%%%%%%%%%%%%%%%%%%%%%%%%%%%%%%%%
 \begin{algorithm}
    \caption{The Ballbot QP-based LPVMPC algorithm}
\begin{algorithmic}[1]
 \renewcommand{\algorithmicrequire}{\textbf{Input:}}
 \renewcommand{\algorithmicensure}{\textbf{Output:}}
 \REQUIRE Initial conditions $x_0$, with reference $x^{\textrm{ref}}_{k},k\in\mathbb{Z}_+\cup\{0\}$
 \ENSURE  The control input $u_k,~k=0,1,\ldots$, that drives the nonlinear system to the reference by satisfying constraints.\\
 \textit{Algorithm steps}
  \STATE Initialize the scheduling parameter for $k=0$ as next:
$$\rho_{i|0}=(\theta_0,\omega_0),~i=0,\ldots,N$$
\vspace{-4mm}
 \WHILE{$k=0,1,\ldots$}
 \STATE Solve the QP in Algorithm \ref{algo:LPVQPOptimProb} 
\begin{equation*}
    u_{i|k}\leftarrow\texttt{QP}(\rho_{i|k},x_k,x_{i+k}^{\text{ref}}),~i=0,\ldots,N-1
    \vspace{-3mm}
\end{equation*}
\STATE Update the scheduling parameters from the designed input $u_{i|k},~i=0,\ldots,N-1$ with  \eqref{eq:dlpv} as $$x_{i+1|k}=\left[\begin{array}{c}
    \phi_{i+1|k} \\
    {\theta_{i+1|k}} \\
    \varphi_{i+1|k} \\
    {\omega_{i+1|k}} 
\end{array}\right]\underset{\text{RK4}}{\overset{\eqref{eq:dlpv}}{=}}f_{\textrm{LPV}}\left(\left[\begin{array}{c}
    \phi_{i|k} \\
    {\theta_{i|k}} \\
    \varphi_{i|k} \\
    {\omega_{i|k}} 
\end{array}\right],u_{i|k}\right)$$
\STATE Update $\rho_{i|k}=\left(\theta_{i|k},\omega_{i|k}\right),~i=0,\ldots,N$
\STATE Apply $u_k=u_{0|k}$ to the continuous system \eqref{eq:lpv} 
\STATE $k\leftarrow k+1$ 
 \STATE Measure $x_{k}$
  \STATE Update $\rho_{i|k}\leftarrow\rho_{i+1|k-1},~i=0,\ldots,N-1$
\ENDWHILE
 \end{algorithmic} 
 \label{algo:lpvmpc}
\end{algorithm}
 %%%%%%%%%%%%%%%%%%%%%%%%%%%%%%%%%%%%%%%%%%%%%%%%%%%%%%%%%%%%%%%%%%%%%%%%%%%%%
As we have introduced in detail the QP minimization problem that is tailored with the LPV formulation in Algorithm \ref{algo:LPVQPOptimProb}, we are ready to proceed by stating the method in Algorithm~\ref{algo:lpvmpc} that is connected with our application on stabilizing and tracking a given reference for the ballbot. Next, every step in the resulting Algorithm~\ref{algo:lpvmpc} is explained. 
 \begin{itemize}
     \item The input to Algorithm \ref{algo:lpvmpc} consists of the initial conditions $x_0$ and the full state reference $x^{\textrm{ref}}$. For the ballbot system, the state reference contains the angular displacement of the ball $\phi_k$. In particular, $\theta_k$ is the tilt angle, and its reference target is zero, thus enforcing stabilization in that operating (unstable equilibrium) point. Finally, for the rest of the states in $x^{\text{ref}}$ that concern the angular velocities, again, we enforce them to zero as the aim is to reach the reference $\phi^{\textrm{ref}}$ and then the ballbot to remain there with $\theta=\dot{\theta}=\dot{\phi}=0$ (steady state) and $\phi$ at the target point.  
     \item In Step 1, initialization at any instant $k$ of the scheduling parameters can be done either from the previous time as $\rho_{i|k}=\sigma(x_{i+1|k-1}),~i=0,\ldots,N-1$ or from the initial conditions as explained in Algorithm \ref{algo:LPVQPOptimProb}. 
     \item In Steps 2-3, we can call the QP Algorithm \ref{algo:LPVQPOptimProb} to provide the optimally designed input $u_{i|k},~i=0,\ldots,N-1$. Note that by having $N$ inputs, we can compute $N+1$ state estimation through the LPV predictor. 
     \item In Step 4, given $x_k$ and $u_{i|k}$, we can compute the true evolution response from the nonlinear system after using the exact/embedded nonlinear discrete predictor \eqref{eq:dlpv}. The scheduling parameters to be used in Step 5 can be revealed in Step 4. The hat notation is omitted as these computations result in the system's true response, assuming that no other disturbance or measurement noise will affect the system.
     \item In Step 6, the one-step implementation of the input to the nonlinear system \eqref{eq:lpv} is done by simulating with ode45 (Runge Kutta 45 with adaptive time step) where between two consequent time instances that differ only one sampling time $t_s$, zero-order hold (ZOH) (constant input of the sampling time) is considered. This high-fidelity simulation ensures good agreement with the true response of the physical plant. 
      \item In Step 7, we slide the prediction horizon by one (sliding the prediction window with length $N$ by one sampling time $t_s$)
     \item In Step 8, we measure the new initial conditions $x_k$.
     \item Finally, in Step 9, we update the new scheduling variable at time $k$ from the previous simulation time $k-1$, and the while loop in Steps 2-10 runs to provide the next control inputs $u_{k}$ for arbitrary $k$.
 \end{itemize}
%%%%%%%%%% RESULTS %%%%%%%%%%%%%%%%%%%%%%%%%%%%%%%%%%%%%%%%%%%%%%%%%%
\section{Results and Discussion}\label{sec:Results}
%%%%%%%%%%%%%%%%%%%%%%%%%%%%%%%%%%%%%%%%%%%%%%%%%%%%%%%%%%%%%%%%%%%%%%%%%%%
\subsection{Refinement of the ballbot physical parameters}\label{sec:Identificationb}
A linearized continuous model at $x_e=[0, 0, 0, 0]^{T}$ can be derived through the LPV evaluated at $x_e$ as in \eqref{eq:linearized}.
\begin{equation}\label{eq:linearized}
    \begin{aligned}
A_c(0,0)& = \begin{bmatrix}
    0 & 0 & 1 & 0 \\
    0 & 0 & 0 & 1 \\
    0 & \underbrace{A_{32}(0,0)}_{p_1} & \underbrace{A_{33}(0,0)}_{p_2} & \underbrace{A_{34}(0,0)}_{0} \\
    0 & \underbrace{A_{42}(0,0)}_{p_4} & \underbrace{A_{43}(0,0)}_{p_5} & \underbrace{A_{44}(0,0)}_{0}
\end{bmatrix},~\\
B_c(0,0)& = \begin{bmatrix}
    0 \\
    0 \\
    \underbrace{B_3(0,0)}_{p_3} \\
    \underbrace{B_4(0,0)}_{p_6}
\end{bmatrix},~d(0)=b_1b_3-(-b_2+\ell r_b)^2.
\end{aligned}
\end{equation}
By introducing the functions $h_i(b):\mathbb{R}^4\rightarrow\mathbb{R},~i=1,~\ldots,6$ with suport the parameterized vector $b=(b_1,b_2,b_3,b_4)$, we can define:
\begin{equation}
\label{eq:eqpb}
    \begin{aligned}
        A_{32}(0,0)&:=h_1(b)=\frac{gl\left(b_{2}-l\mathrm{r_b}\right)}{d(0)},\\
        A_{33}(0,0)&:=h_2(b)=-\frac{b_{3}\,b_{4}}{d(0)},\\
        A_{42}(0,0)&:=h_4(b)=\frac{b_{1}\,g\,l}{d(0)},\\
        A_{43}(0,0)&:=h_5(b)=-\frac{b_{4}\,\left(b_{2}-l\,\mathrm{r_b}\right)}{d(0)},\\
        A_{34}(0,0)&:=0,~A_{44}(0,0):=0,\\ 
        B_3(0,0)&:=h_3(b)=\frac{b_{3}\,\mathrm{r_b}-\mathrm{r_b}\,\left(b_{2}-l\,\mathrm{r_b}\right)}{\mathrm{r_w}\,d(0)},\\
        B_4(0,0)&:=h_6(b)=\frac{\mathrm{r_b}\,\left(b_{2}-l\,\mathrm{r_b}\right)-b_1\mathrm{r_b}}{\mathrm{r_w}\,d(0)}.
    \end{aligned}
\end{equation}
To infer the nonlinear $b$-parameters, given the $p$-parameters values, we must solve the nonlinear system of equations described in the set of equations \eqref{eq:eqpb} for the unknown vector $b=(b_1,b_2,b_3,b_4)$ and with the functionality $h_i(b):\mathbb{R}^4\rightarrow\mathbb{R},~i=1,\ldots,6$. The problem remains to find $b^*$ that satisfies $p_i-h_i(b^*)\approx 0$ for all $i=1,\ldots,6$. Thus, we define the operator $F(b):\mathbb{R}^4\rightarrow\mathbb{R}^6$ along with the Jacobian next:
\begin{equation*}
    \begin{aligned}
        F(b):=\left[\begin{array}{cc}
             p_1-h_1(b)  \\[1mm]
              p_2-h_2(b)  \\[1mm]
              \vdots\\[1mm]
              p_6-h_6(b)  \\
        \end{array}\right],~\nabla F(b)=\left[\begin{array}{ccc}
            \frac{\partial F_1}{\partial b_1} & \cdots & \frac{\partial F_1}{\partial b_4}\\[1mm]
           \frac{\partial F_2}{\partial b_1} & \cdots & \frac{\partial F_2}{\partial b_4}\\[1mm]
           \vdots\\[1mm]
           \frac{\partial F_6}{\partial b_1} & \cdots & \frac{\partial F_6}{\partial b_4}\\
        \end{array}\right].
    \end{aligned}
\end{equation*}
The Newton scheme for updating the parameters $b$ with initialization $b_0$ and for $k\in\mathbb{Z}_+\cup\{0\}$ is 
\begin{equation}
\label{eq:Newton}
    b_{k+1}=b_{k}-\left[\nabla F(b_k)\right]^{-1}F(b_k),~k=0,\ldots
\end{equation}
Upon convergence of \eqref{eq:Newton}, the optimal solution vector $b^*$ can be obtained to satisfy $F(b^*)\approx 0$ (with the smallest possible residual). In such problems, where the nonlinear parameterized manifold is nonconvex, the iterative scheme in \eqref{eq:Newton}  remains divergent when random initialization is considered, or a solution can be algebraically correct without explaining the physical law. To tackle such a problem, initialization of the $b$-parameters should be done by respecting the engineering regularity of the problem. Thus, some parameters $b_i$ can be initialized with physical meaning, e.g., $b_1,~b_2,~b_3$ being positive to explain the high-fidelity knowledge of the physical plant \cite{Max} where $b_4$ as concerns friction can remain arbitrary but also with a reasonable small value. 

\begin{figure}[!h]
    \centering
    \includegraphics[scale=0.18]{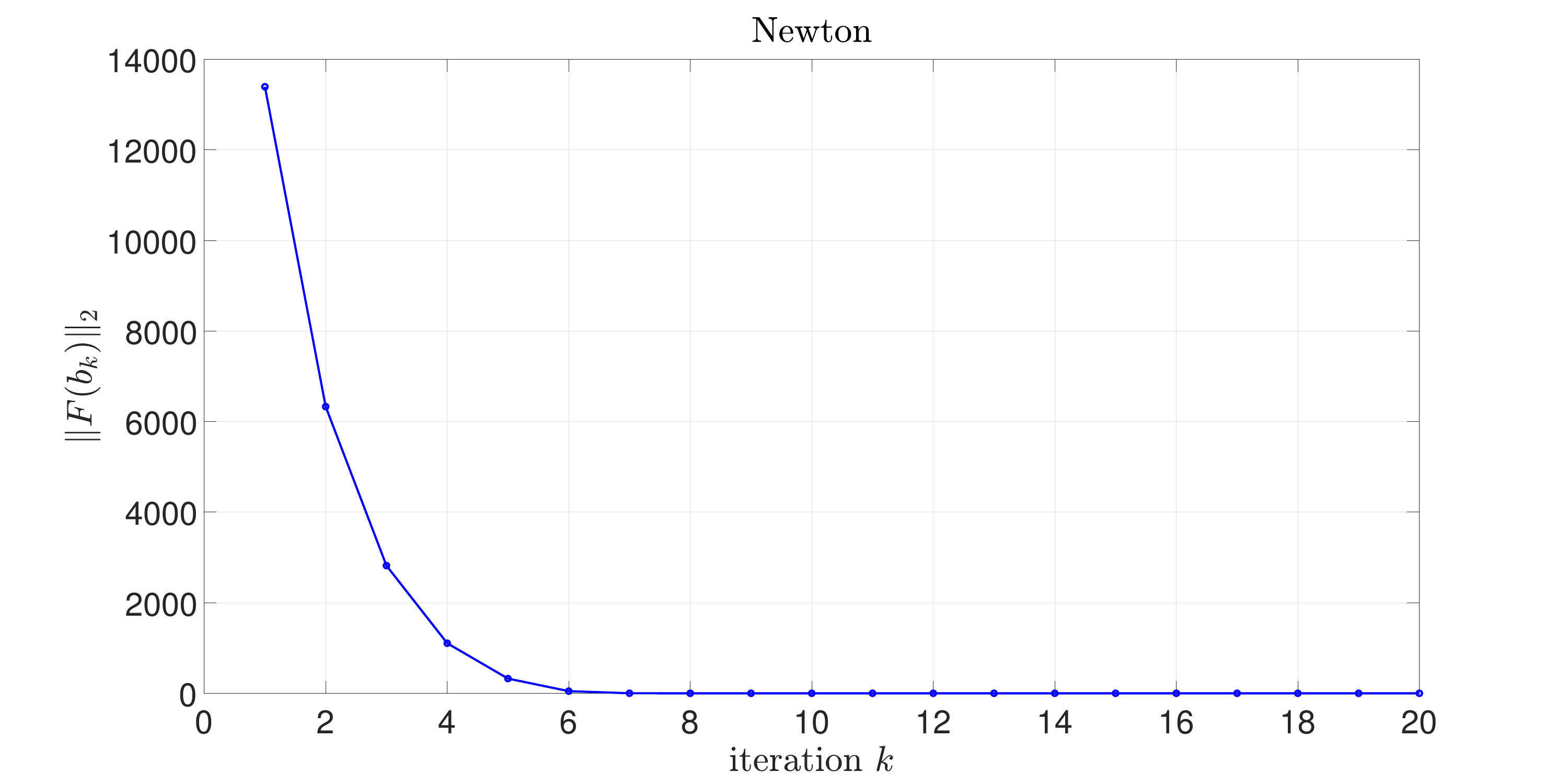}
    \caption{Convergence of the Newton scheme with residual error $\Vert F(b^*)\Vert=0.4330$. The optimal parameter vector $b^*$ is shown in Tab.~\ref{tab:bparameters}.}
    \label{fig:Newton}
\end{figure}

In Fig.~\ref{fig:Newton}, the iterative scheme \eqref{eq:Newton} converged after five steps and gave the optimal solution with a residual error that stagnated to the value $0.4330$. This error cannot be improved further as the analysis relies on the identified parameters $p$ from the study \cite{Max}, where the dynamics are corrupted with noise and/or slight nonlinear behavior from the physical plant. 

\begin{figure}[!h]
\centering
	\includegraphics[scale=0.18]{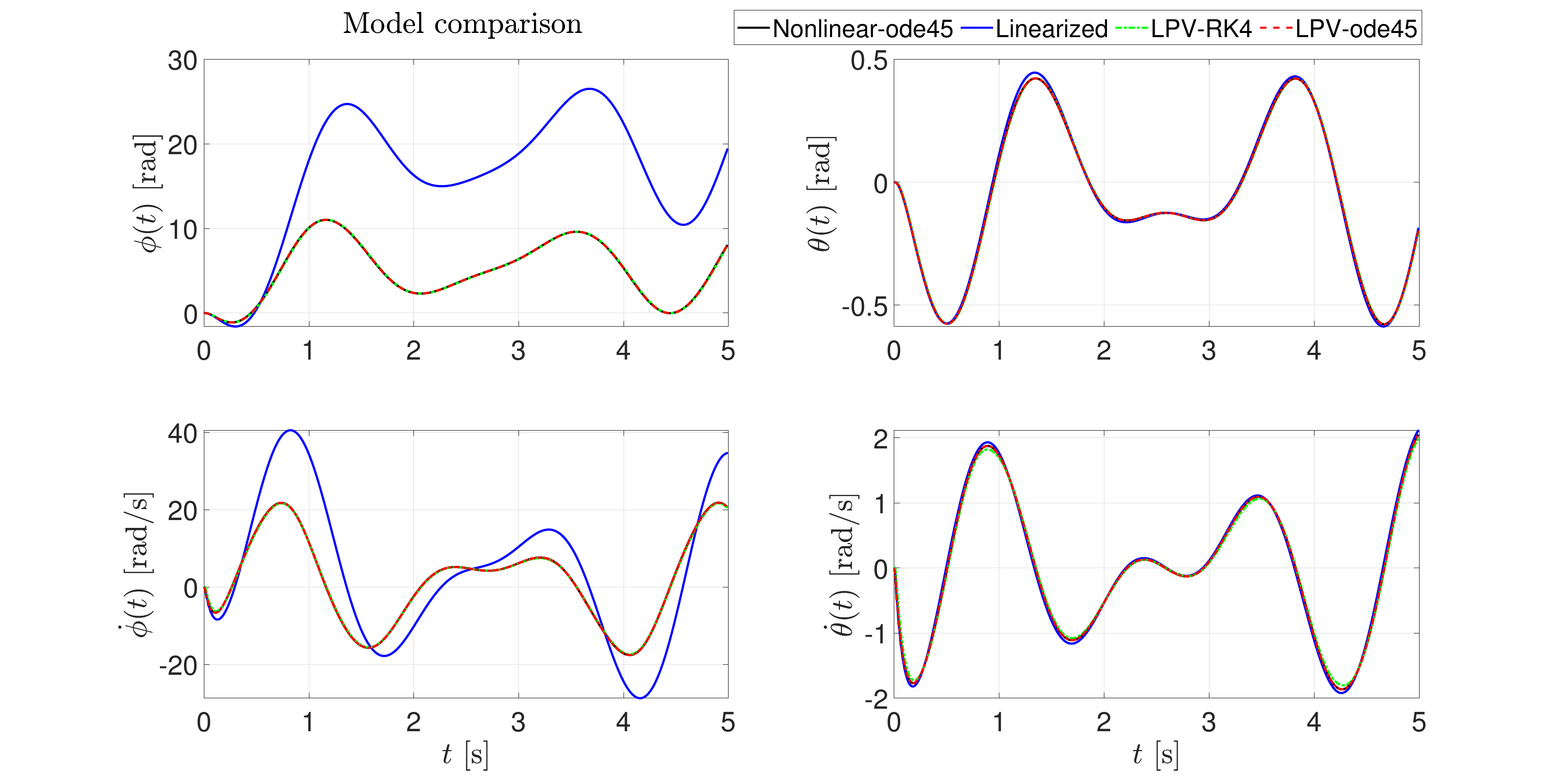}
		\caption{Models simulations; the linearized model given in dashed blue line; non-linear model simulation given in black; LPV in dashed red lines; all done using ode45 on MATLAB. Finally, with green is the RK4 method utilizing \eqref{eq:dlpv}. The input $\tau_y$ is a multiharmonic signal.}
        \label{fig:sim}
\end{figure}

Continuously, in Fig.~\ref{fig:sim}, the simulation indicates the expected performance after substituting the obtained parameters $b$ in Table~\ref{tab:bparameters} and applying a multi-harmonic input to the stabilized plant. In particular, in Fig.~\ref{fig:sim}, the nonlinear model compared to the linearized model stays close to the dynamical evolution for all states when the dynamics are close to the linearization operation point $x_e=0$ as expected. The discrepancy increases for larger deviations in $\phi$ from the origin where the linearized cannot be trusted. Finally, in Fig.~\ref{fig:sim}, the comparison between the LPV and the nonlinear model certifies that the LPV is an equivalent embedding of the nonlinear system where differences are buried in machine precision error. Finally, the discrete-time implementation with the RK4 in \eqref{eq:dlpv} with a ZOH also certifies a good accuracy compared to the ode45.

\begin{table}[!h]
\centering
\caption{The refined $b$ parameters of the nonlinear model.}
\footnotesize
\begin{tabular}{|c |c |c |c | c |}
\hline
 $b^*$ & $b_1$ & $b_2$ & $b_3$ & $b_4$   \\ [1ex] 
 Optimal & 0.002483 & 0.059325 & 0.143093 & -0.07436  \\  [1ex] 
 \hline 
\end{tabular}
\label{tab:bparameters}
\end{table}

\subsection{LPVMPC for Reference Tracking}\label{sec:reference_tracking}
Reference tracking for a ballbot refers to the stabilized motion by the LPVMPC method that can maneuver over any given state reference $x^{\mathrm{ref}}$. We will consider three study scenarios. The 1st scenario will be a 1D reference tracking of a nonsmooth discontinuous trajectory that will "shock" the optimization problem due to the sudden changes. The 2nd scenario will be a 2D smooth reference as Lissajous curves, assuming the ballbot can be driven independently in the $xz$ and $yz$ planes. Finally, a 3rd scenario will illustrate the reference tracking of a single set point with a scheme of the LPVMPC that can guarantee the stability of the closed-loop system and the recursive feasibility of its optimization problem. These properties are explained below. Tuning of the prediction horizon depends on the physical system's operational bandwidth and can vary significantly across applications. The ballbot is a robotic system that operates around $1$(Hz); thus, a horizon of about $1$ (s) is reasonable.
\subsubsection{1st Scenario: Nonsmooth 1D reference}
\begin{figure}[!h] 
\centering
   \includegraphics[scale=0.18]{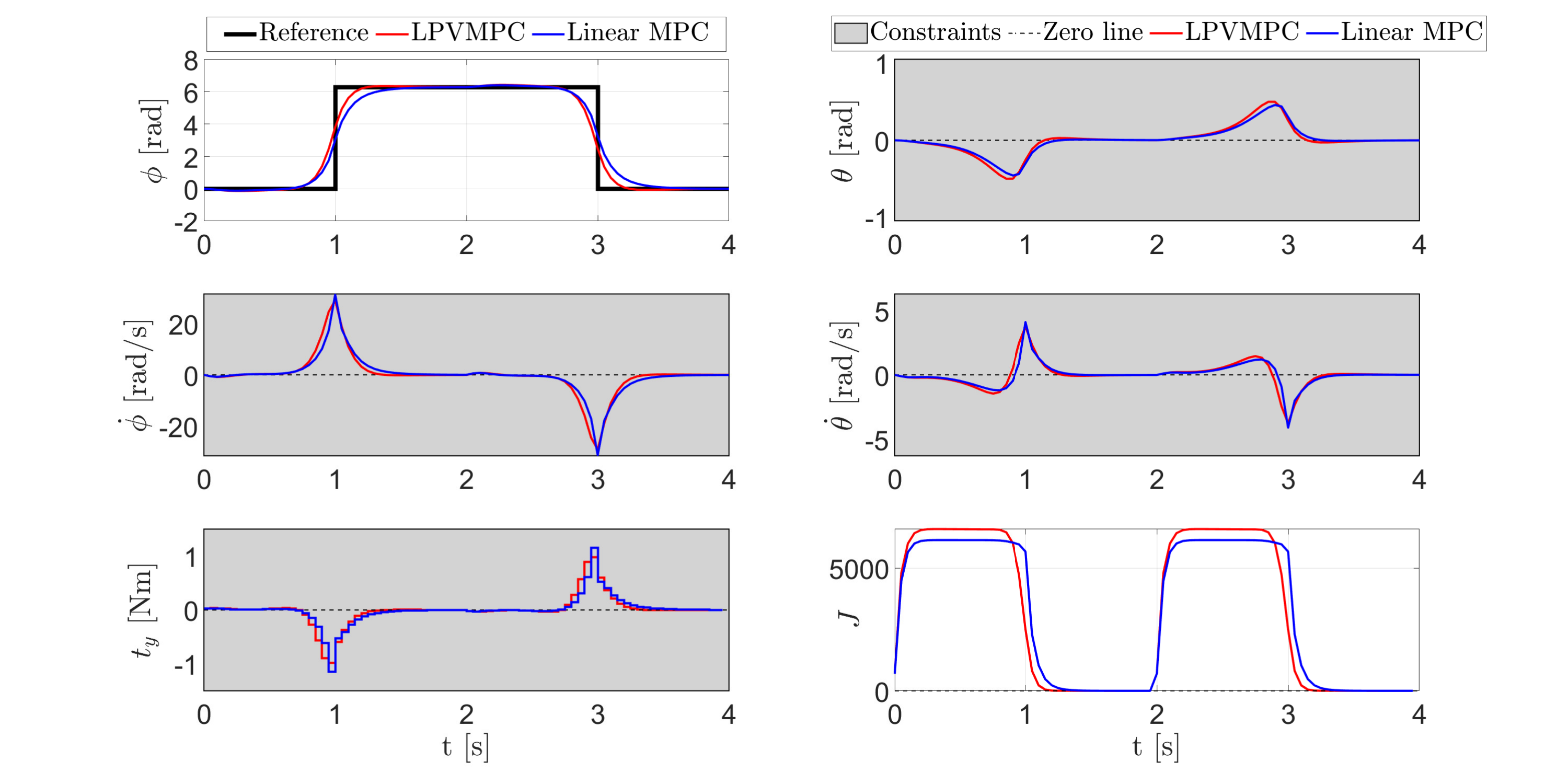}
   \caption{Ballbot reference tracking with the angular displacement $\phi$ traveling from the origin to $2\pi$ rad and back. A comparison between the linear and LPV MPC frameworks is illustrated. $\theta$ measuring balancing of ballbot, $\dot{\phi}$, $\dot{\theta}$ measuring the angular velocities, $\tau_y$ is the virtual torque, $J$ is the energy cost function.}
   \label{fig:1Dref} 
\end{figure}

In Fig.~\ref{fig:1Dref}, we consider the following 2-set point nonsmooth reference. The ballbot should start from the origin and at $t=1$ s should roll for one circle $2\pi$ rad, then stay there for $2$ s and for the remaining $1$ s to return to the origin; thus, it will roll a total distance of approximately $(2\cdot 2\pi\cdot r_b\approx 1.5~\text{m})$ within $4$ s. 

At the beginning and near the origin, positive virtual torque applies on the ballbot to decrease the tilt angle $\theta$ that instantly will move the ballbot to the negative $\phi$, where immediately switching direction (normal behavior of such non-minimal phase system, similar to a bike where to turn, first a small maneuver to the opposite direction is needed). Continuously the ballbot increases speed in the $\phi$ positive direction while approaching the target $2\pi$ with satisfying accuracy and faster than the linear one after $\approx 1.3$ s, then stays there and at around $2^+$ s where the prediction horizon of time length $N\cdot t_s=20\cdot 0.05=1$ s receives the information of the sudden change in reference from $2\pi$ to $0$ rad, the ballbot slightly overshoots in $\phi$ so to change the tilt angle again and starts to accelerate in the opposite direction. After less than a second, the ballbot has reached its origin without overshooting, which outlines the good performance seen mainly in nonlinear MPC frameworks. 

In Fig.~\ref{fig:1Dref}, we compare the MPC performances between the LPV and the linearized model. As expected, the LPVMPC is generally faster at reaching the reference with almost similar computational complexity as the linear MPC and the potential to handle better strong maneuvering as the adaptive scheduling variable $\rho_{i|k}$ provides linearizations over a trajectory of the nonlinear ballbot model instead of a single fixed point for the linearized model. The quadratic costs for the 1st scenario are defined as $Q = \texttt{diag}([200,1,0.1,0.1])$, $R = 1000$, and the terminal cost $P$ is the Riccati solution from the LQR on the linearized model from \cite{Max} with the same tuning of $Q$ and $R$ that will penalize the "tail" of the horizon. The sampling time is $t_s = 0.05$ (s), with horizon length $N=20$, and the MPC solution is obtained under the input/state hard constraints introduced in Table~\ref{tab:MPCconstraints} and highlighted with the gray background in Fig.~\ref{fig:1Dref}. 
\begin{table}[!h]
    \centering
    \caption{MPC Constraints}
    \label{tab:MPCconstraints}
    \setlength{\tabcolsep}{2pt}
    \begin{tabular}{c | c c |c} 
       \textbf{Variable}   & \textbf{Lower bound}  &  \textbf{Upper bound} & \textbf{Units}\\ \hline
                         $\phi_k$ &   $-\infty$ &       $+\infty$ & [rad]  \\
        $\theta_k$ &   $-\pi/3$  &      $+\pi/3$ &[rad] \\
          $\dot{\phi}_k$ &   $-10\pi$  &     $+10\pi$ & $[\frac{\rm rad}{\rm s}]$  \\
           $\dot{\theta}_k$ &   $-2\pi$  &    $+2\pi$ & $[\frac{\rm rad}{\rm s}]$ \\\hline
              $u$ &   $-1.5$ &    $+1.5$ & [\text{Nm}]   
    \end{tabular}
\end{table}
In the 1D reference tracking, as shown in Fig.~\ref{fig:1Dref}, traveling in $y$ direction, the ballbot system also balances itself simultaneously, as shown with $\theta$ being between ranges of $0 \,rad$ to $0.5\,rad$ while traveling a distance of $\phi$, which is referenced for $2\pi \, rad$, and angular velocities of both are shown in $\dot{\theta}$, and $\dot{\phi}$ respectively. They show the agile and robust movements of the ballbot. Also, in correspondence to the reference and torque, $\tau_y$ is applied to the virtual wheel actuator. The virtual torque is shown as discrete steps, i.e., ZOH, and the energy function $J$ to be minimized to get the optimal control input $\tau_{y}$.
\subsubsection{2nd Scenario: Lissajous 2D smooth curves}
\begin{figure}[!h] 
\centering
   \includegraphics[scale=0.18]{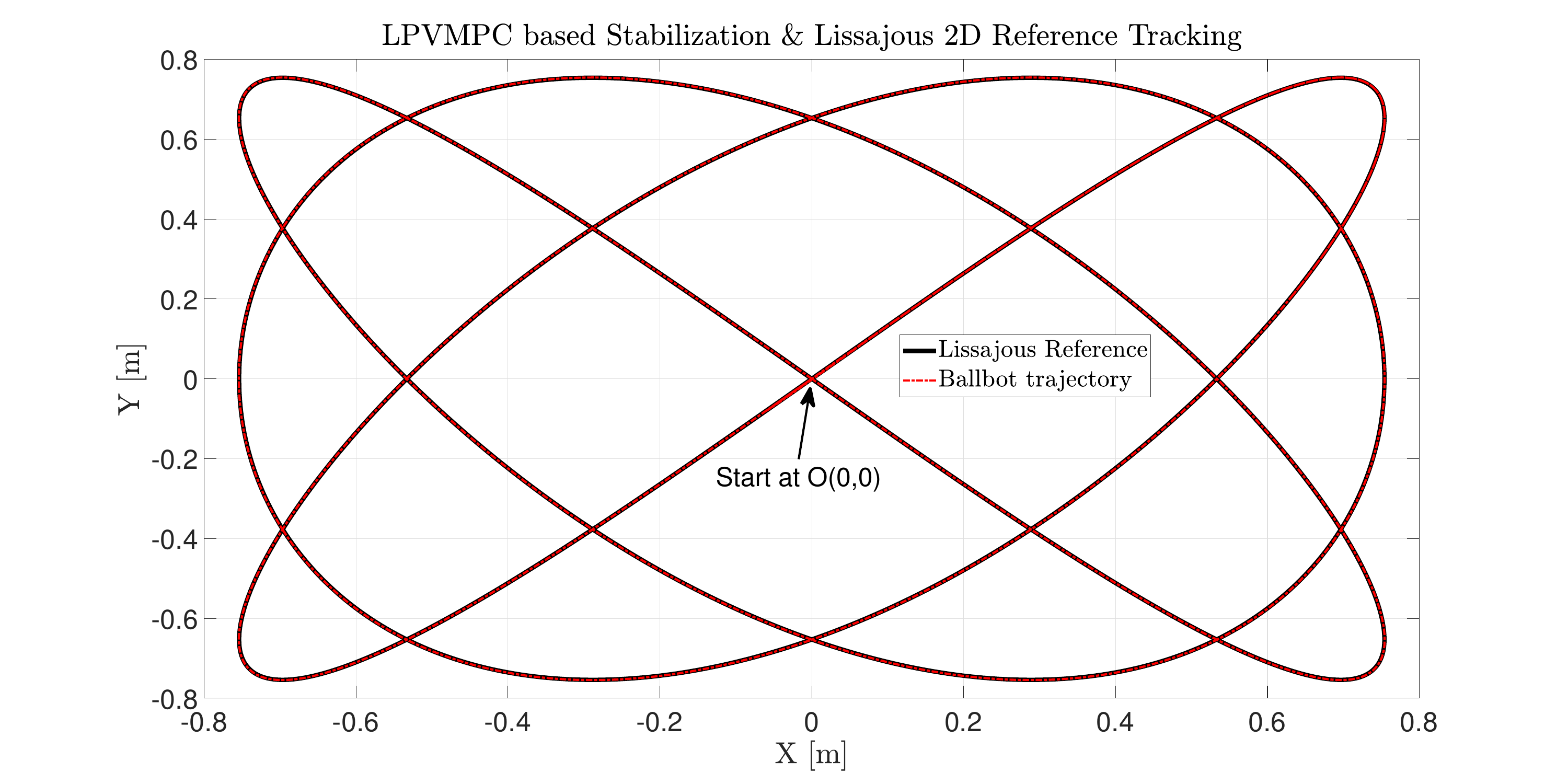}
   \caption{Stabilization and Lissajous curve trajectory tracking of the ballbot in 2D using the LPVMPC algorithm solved in both $xz$ and $yz$ planes.}
   \label{fig:mpc1} 
\end{figure}

Solving the LPVMPC reference tracking problem in parallel for the $yz$-plane and $xz$-plane, we can drive the ballbot over the $xy$ plane with reference as Lissajous curves in Fig.~\ref{fig:mpc1}. The ballbot starts from $(x,y) = (\phi_{x}r_b,\phi_{y}r_b) = (0,0)$, and following the path to global coordinates as: $X^{\textrm{ref}}=\phi_x r_b = 0.12\cdot 2\pi\sin(0.3t)$ and $Y^{\textrm{ref}}=\phi_y r_b=0.12\cdot 2\pi\sin(0.4t)$, under the constraints in Table \ref{tab:MPCconstraints}. In that case, the reference is smoother without sudden changes. Therefore, we can increase the weight for matching the $\phi$ without inheriting infeasibility problems; thus, the quadratic cost for the 2D reference tracking in each direction can be set as $Q = \texttt{diag}([1000,1,0.1,0.1])$. All the other quadratic weights and parameter specifications are considered the same as in the 1st scenario. In Fig.~\ref{fig:mpc1}, the LPVMPC control strategy accurately drives the ballbot to the reference within an area of around $2~m^2$ and for $70$ s. 
\subsubsection{3rd Scenario: Single set point reference and LPVMPC with guarantees}
In our last scenario, regulated as in the 1st (e.g., quadratic costs, etc.), we want to illustrate the crucial advantage of considering LPV models for MPC tasks, as it can provide theoretical guarantees. On the one hand, when a linearized model is considered instead of the underlying nonlinear, achieving theoretical guarantees within MPC could be a straightforward task. Still, it will remain quite conservative and thus not practically useful. On the other hand, in the case of NMPC, guarantees can also be provided at the expense of computational complexity. Consequently, an excellent alternative for providing theoretical guarantees is the LPV formulation within MPC, which can be considered an excellent trade-off between the complexity and conservativeness of the previous two general approaches. 

\begin{figure}[!h]
    \centering
    \includegraphics[scale=0.18]{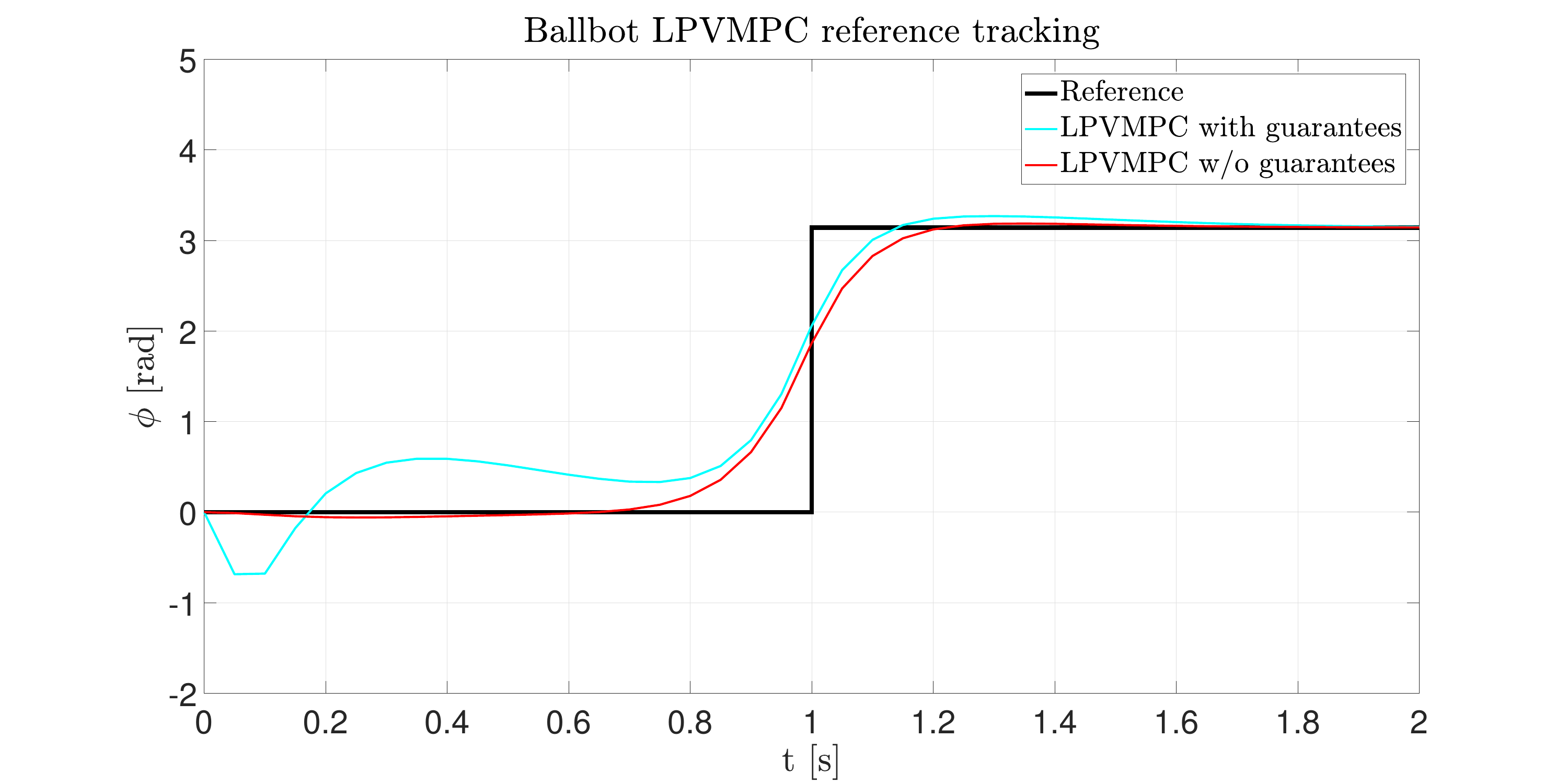}
    \caption{Single-set point reference i.e., $\pi$, for the LPVMPC framework with and without theoretical guarantees.}
    \label{fig:guarantees}
\end{figure}

To prove theoretical guarantees, in Algorithm \ref{algo:LPVQPOptimProb}, we have to introduce additionally the terminal equality constrain as 
\begin{equation}
\label{eq:terminalconstraint}
    \mathcal{X}_{t}=\{\hat{x}_{N|k}\in\mathbb{R}^{n_{\textrm{x}}}~|~\hat{x}_{N|k}=x^{\text{ref}}_{N+k}\}.
\end{equation}
Thus, the state constraints $\mathcal{X}$ in Algorithm \ref{algo:LPVQPOptimProb} are enhanced with $\hat{x}_{N|k}\in\mathcal{X}_t$. Note that in the presence of the terminal constraint $\mathcal{X}_t$, the terminal cost remains zero; thus, in that case, no use of the weight $P$ is needed. In Fig.~\ref{fig:guarantees}, the dynamics are illustrated when considering the terminal constraint (cyan line) in \eqref{eq:terminalconstraint}, where we lose some performance compared with no guarantees (red line). Thus, with active \eqref{eq:terminalconstraint} and making use of the Theorem 1 in \cite{verhoek2023linear} (page: 6), all the conditions are satisfied in our 3rd scenario, which leads to the following conclusions;
i) the LPVMPC optimization problem is recursively feasible; ii) the closed-loop system satisfies the constraints $u_{i|k}\in\mathcal{U}$, $\hat{x}_{i|k}\in\mathcal{X}$, and $\hat{x}_{N|k}\in\mathcal{X}_t$; iii) $(u^\text{ref},x^\text{ref})=(0,\pi)$ is an exponentially stable (forced) equilibrium of the closed-loop system.
The proof of these properties follows the same reasoning as the proof of Theorem 1 in \cite{verhoek2023linear}.
\begin{table}[!h]
    \centering
    \caption{Comparison of computation times for a single QP using YALMIP with LPVMPC. The simulations are performed on a Dell Latitude 5590 laptop with an Intel(R) Core(TM) i7-8650U CPU and 16 GB of RAM. The scenarios are implemented in MATLAB, utilizing the YALMIP toolbox~\cite{Lofberg2004}, with an optimality tolerance of $10^{-8}$.}
    \label{tab:Compuation_time}
    \setlength{\tabcolsep}{2pt}
    \begin{tabular}{ c | c | c} 
        \toprule
        \rowcolor{mydarkgray} \textbf{LPVMPC with $t_s=0.05$ s}   & \textbf{Mean}  & \textbf{Standard deviation}  \\ \hline
                             1st Scenario &  0.0086 s &  0.0104 s \\
        \rowcolor{mygray}    2nd Scenario &  0.0087 s  & 0.0132 s \\
                             3rd Scenario & 0.0058 s  & 0.0025 s   
    \end{tabular}
\end{table}
%%%%%%%%%%%%%%%%%%%%%%%%%%%%%%%%%%%%%%%%%%%%%%%%%%%%%%%%%%%%%%%%%%%%%%%%%%%%%%%%%%%
% CONCLUSION
\section{Conclusion}
%%%%%%%%%%%%%%%%%%%%%%%%%%%%%%%%%%%%%%%%%%%%%%%%%%%%%%%%%%%%%%%%%%%%%%%%%%%%%%%%%%%
In this study, we were concerned with the reference tracking of a ballbot by utilizing the LPV embedding within the MPC framework. An important step was the refinement of the parameters of the ballbot system from measurements of a past study. Thus, having the nonlinear model allowed the nonlinear embedding of the LPV form. An improved discretization method based on the well-known Runge-Kutta $4^{\textrm{th}}$ order introduced in the LPVMPC framework potentially improves the accuracy issue that can arise in the forward Euler method. The results indicated that applying the LPVMPC-based control method for simultaneously balancing and reference tracking can achieve real-time implementation as the average timing of solving the QP problem is usually $(10)$-times less than the sampling time under consideration.

Utilizing such a control design that preserves the nonlinear behavior of the model is advantageous and practical, unlike methods that work with linearized models or try to solve very complex optimization problems tailored directly to the nonlinear, which might be impossible online. The method is illustrated with practical reference tracking scenarios where theoretical guarantees can be provided in the single-set point case. Proving stability and recursive feasibility in the case of multiple reference points is left to future investigation. Finally, as the theoretical establishment of such a robotic system has matured, implementing this method in the physical ballbot in a real environment is left for our immediate future endeavors.
%%%%%%%%%%%%%%%%%%%%%%%%%%%%%%%%%%%%%%%%%%%%%%%%%%%%%%%%%%%%%%%%%%%%%%%%%%%%%%%%%%%
% ACKNOWLEDGEMENT
\section*{Acknowledgement}
We thank Mr. Ing.-Ievgen Zhavzharov, who constructed the ballbot in Fig.~\ref{fig:ballbot_SideView}.
%%%%%%%%%%%%%%%%%%%%%%%%%%%%%%%%%%%%%%%%%%%%%%%%%%%%%%%%%%%%%%%%%%%%%%%%%%%%%%%%%%%
% AUTHORS’ STATEMENT
%\section*{Authors' Statement}
%%%%%%%%%%%%%%%%%%%%%%%%%%%%%%%%%%%%%%%%%%%%%%%%%%%%%%%%%%%%%%%%%%%%%%%%%%%%%%%%%%%
% BIBLIOGRAPHY: Write list of literature in literature.bib in bibstyleformat (http://www.bibtex.org/)
\bibliographystyle{IEEEtran}
\bibliography{main}
%%%%%%%%%%%%%%%%%%%%%%%%%%%%%%%%%%%%%%%%%%%%%%%%%%%%%%%%%%%%%%%%%%%%%%%%%%%%%%%%%%%
\end{document}